\documentclass[preprint,times]{elsarticle}

\journal{Linear Algebra and its Applications}

\usepackage{amsmath}
\usepackage{amscd,amssymb}

\begin{document}

\begin{frontmatter}

 \title{On powers of Stirling matrices}
 \author{Istv\'an Mez\H{o}}
 \address{Department of Algebra and Number Theory, Institute of Mathematics, University of Debrecen, Hungary\\University of Debrecen, H-4010, Debrecen, P.O. Box 12, Hungary}
 \ead{imezo@math.klte.hu}
 \ead[url]{http://www.math.klte.hu/algebra/mezo.htm}

\begin{abstract}
The powers of matrices with Stirling number-coefficients are investigated. It is revealed that the elements of these matrices have a number of properties of the ordinary Stirling numbers. Moreover, ''higher order" Bell, Fubini and Eulerian numbers can be defined. Hence we give a new interpretation for E. T. Bell's iterated exponential integers. In addition, it is worth to note that these numbers appear in combinatorial physics, in the problem of the normal ordering of quantum field theoretical operators.
\end{abstract}

\end{frontmatter}

\newcommand{\stirlingf}[2]{\genfrac[]{0pt}{}{#1}{#2}}
\newcommand{\stirlings}[2]{\genfrac\{\}{0pt}{}{#1}{#2}}

\newcommand{\ZZ}{\mathbb{Z}}
\newcommand{\RR}{\mathbb{R}}
\newcommand{\mS}{\mathcal{S}}
\newcommand{\im}{\mbox{\upshape Im}}

\newtheorem{Theorem}{Theorem}
\newtheorem{Lemma}[Theorem]{Lemma}
\newtheorem{Definition}[Theorem]{Definition}
\newtheorem{Proposition}[Theorem]{Proposition}
\newtheorem{Corollary}[Theorem]{Corollary}

Stirling cycle numbers (also known as Stirling number of the first kind \cite{GKP}) count the number of permutations with given cycles. Concretely, $|s_{n,m}|$ (also denoted as $\stirlingf{n}{m}$ or $c(n,m)$) gives the number of $n$-paritions with exactly $m$ cyc\-les $(0\le m\le n)$. The signed Stirling cycle number is $(-1)^{n-m}|s_{n,m}|=:s_{n,m}$. The exponential generating function for the latter is
\[\sum_{n=0}^\infty s_{n,m}\frac{x^n}{n!}=\frac{1}{m!}(\ln(1+x))^m.\]

Stirling set number (also known as Stirling number of the second kind \cite{GKP}), denoted by $S_{n,m}$ (or $\stirlings{n}{m}$), enumerates the partitions of an $n$-set into $m$ subsets $(0\le m\le n)$. The exponential generating function for them is
\[\sum_{n=0}^\infty S_{n,m}\frac{x^n}{n!}=\frac{1}{m!}\left(e^x-1\right)^m.\]

We enumerate the main identities involving these numbers. Let $(x)_n=x(x-1)\cdots(x-n+1)$ (with $(x)_0=1$), then the next matrix identities are valid \cite{GKP}
\begin{equation}
\begin{pmatrix}s_{0,0}&0&0&\cdots&0\\s_{1,0}&s_{1,1}&0&\cdots&0\\\vdots&&\ddots\\s_{n,0}&s_{n,1}&\cdots&s_{n,n-1}&s_{n,n}\end{pmatrix}\begin{pmatrix}x^0\\x^1\\\vdots\\x^n\end{pmatrix}=\begin{pmatrix}(x)_0\\(x)_1\\\vdots\\(x)_n\end{pmatrix},\label{mat1}
\end{equation}
and
\begin{equation}
\begin{pmatrix}S_{0,0}&0&0&\cdots&0\\S_{1,0}&S_{1,1}&0&\cdots&0\\\vdots&&\ddots\\S_{n,0}&S_{n,1}&\cdots&S_{n,n-1}&S_{n,n}\end{pmatrix}\begin{pmatrix}(x)_0\\(x)_1\\\vdots\\(x)_n\end{pmatrix}=\begin{pmatrix}x^0\\x^1\\\vdots\\x^n\end{pmatrix}.\label{mat2}
\end{equation}
So the Stirling matrices are inverses of each other.

Applying again the matrix on \eqref{mat2}, we get the Bell polynomials \cite{Aigner}:
\begin{equation}
\begin{pmatrix}S_{0,0}&0&0&\cdots&0\\S_{1,0}&S_{1,1}&0&\cdots&0\\\vdots&&\ddots\\S_{n,0}&S_{n,1}&\cdots&S_{n,n-1}&S_{n,n}\end{pmatrix}\begin{pmatrix}x^0\\x^1\\\vdots\\x^n\end{pmatrix}=\begin{pmatrix}B_0(x)\\B_1(x)\\\vdots\\B_n(x)\end{pmatrix}.\label{Stoxpow}
\end{equation}
That is,
\begin{equation}
B_n(x)=\sum_{m=0}^n S_{n,m}x^m.\label{Bellsum}
\end{equation}
As one can see by the definition, $B_n:=B_n(1)$ gives all the possible partitions of a set with $n$ element. These numbers are called Bell numbers \cite{Aigner}. Among others, we point out that powers of the Stirling matrices present new numbers which satisfies all the relevant properties of the Bell numbers. It seems that these numbers were not investigated except by E. T. Bell himself \cite{Bell}.

In all what follows, we would like to discuss the Stirling matrices on a unified base as far as possible. So we introduce the notation
\begin{equation}
S^{-1}_{n,m}:=s_{n,m},\label{Ssnot}
\end{equation}
and the infinite lower triangular matrices
\[\mS=\begin{pmatrix}S_{0,0}&0&0&0&0&\cdots\\S_{1,0}&S_{1,1}&0&0&0&\cdots\\S_{2,0}&S_{2,1}&S_{2,2}&0&0&\cdots\\\vdots&\vdots&\ddots&\vdots&\vdots\\S_{n,0}&S_{n,1}&\cdots&S_{n,n}&0&\cdots\\\vdots\end{pmatrix},\]
and $\mS^{-1}$, which is the inverse of $\mS$. That is, $(\mS^{-1})_{n,m}=S^{-1}_{n,m}$, with the notation \eqref{Ssnot}.

The formal $p$-th power of these matrices will be denoted as $\mS^p$ and the inverse is $\mS^{-p}$, while their $(n,m)$-th entries are $S^p_{n,m}$ and $S^{-p}_{n,m}$, respectively.

\section{Basic properties of the entries of $\mS^{\pm p}$}

First, to discuss the generating functions on the most simple way, we introduce the functions 
\[\sigma^1(x):=e^x-1,\quad\sigma^p(x):=\sigma^{p-1}(\sigma^{1}(x)),\]
and
\[\sigma^{-1}(x):=\ln(1+x),\quad\sigma^{-p}(x):=\sigma^{-(p-1)}(\sigma^{-1}(x)),\]
here $p>1$. Further, let $\sigma^0(x)=x$.

Then we have the general
\begin{Proposition}\label{gfs}The generating function of the $m$-th column of $\mS^p$ ($p\in\ZZ$) is
\[\sum_{n=0}^\infty S^p_{n,m}\frac{x^n}{n!}=\frac{(\sigma^p(x))^m}{m!}.\]
\end{Proposition}

\textit{Proof.} For $p\in\{-1,0,1\}$ this is known or trivial. Let $p>1$. Then
\[\sum_{n=0}^\infty S^p_{n,m}\frac{x^n}{n!}=\sum_{n=0}^\infty \left(\sum_{k=0}^nS^1_{n,k}S^{p-1}_{k,m}\right)\frac{x^n}{n!}=\sum_{k=0}^\infty S^{p-1}_{k,m}\left(\sum_{n=0}^\infty S^{1}_{n,k}\frac{x^n}{n!}\right)=\]
\[\sum_{k=0}^\infty S^{p-1}_{k,m}\frac{(\sigma^{1}(x))^k}{k!}=\frac{[\sigma^{p-1}(\sigma^{1}(x))]^m}{m!},\]
by induction. The case of $p<-1$ is totally similar.
\hfill\qed

\section{Higher order Bell numbers}

Next we show that the well known identity \cite{GKP}
\begin{equation}
S_{n,m}=\frac{1}{m!}\sum_{k=0}^m\binom{m}{k}(-1)^{m-k}k^n\label{Skadn}
\end{equation}
can be generalized and thus the ``higher order" Bell numbers and polynomials defined below appear naturally. To see this, we cite Bell's definition \cite{Bell} for the so-called ``iterated exponential integers" -- mutatis mutandis.

\begin{Definition}[E. T. Bell -- 1938]The higher order Bell polynomials (or iterated exponential polynomials) can be defined via their exponential generating functions. Let
\[E^{p}(t):=e^{\sigma^p(t)},\]
where $p\ge0$. Then we can define the polynomials
\begin{equation}
B_n^p(x):=n![t^n](E^p(t))^x,\label{Belldef}
\end{equation}
and the ``higher order" Bell numbers
\[B_n^p:=B_n^p(1).\]
\end{Definition}
It can be seen that $B_n^0(x)=x^n$.

We give the following table of the numbers $B_n^p$. In the paper of Ginsburg \cite{Gins} there is a misprint. $B_7^3$ was given incorrectly by him.
\begin{center}
\begin{tabular}{r|cccccccc}$B_n^p$&$n=1$&$n=2$&$n=3$&$n=4$&$n=5$&$n=6$&$n=7$&\\\hline
$p=1$&1&2&5&15&52&203&877\\
$p=2$&1&3&12&60&358&2471&19 302\\
$p=3$&1&4&22&154&1304&12 915&146 115\\
$p=4$&1&5&35&315&3455&44 590&660 665\\
$p=5$&1&6&51&561&7556&120 196&2 201 856\\
\end{tabular}
\end{center}

\vspace{.2cm}

We remark that Bell used the notation and indices $\xi_n^{(m)}$ instead of $B_n^p$ but the latter one fits in with our approach. He did not use the polynomials but only the numbers and gave several properties for $B_n^p$. On the other hand, Ginsburg used the letter $E$ also for the generating functions of the iterated exponentials but he applied lower indices.

It should also be noted that these numbers appear in a paper concerning combinatorial physics: normal ordering of boson annihilation and creation operators \cite{Blasetal}.

With these considerations, we present the next

\begin{Proposition}\label{PropBellSt}For $p\ge0$ we have the more general form of \eqref{Stoxpow} as
\begin{equation}
\begin{pmatrix}B_0^p(x)\\B_1^p(x)\\\vdots\end{pmatrix}=\mS^p\begin{pmatrix}x^0\\x^1\\\vdots\end{pmatrix}.\label{BellStsum}
\end{equation}
and for $p\ge1$ we get the generalization of \eqref{Skadn}:
\[S_{n,m}^p=\frac{1}{m!}\sum_{k=0}^m\binom{m}{k}(-1)^{m-k}B_n^{p-1}(k).\]
\end{Proposition}

\textit{Proof.} We consider the series
\[\sum_{n=0}^\infty\sum_{m=0}^\infty S_{n,m}^{p+1}\cdot(x)_m\frac{t^n}{n!}=\sum_{m=0}^\infty(x)_m\sum_{n=0}^\infty S_{n,m}^{p+1}\frac{t^n}{n!}=\sum_{m=0}^\infty(x)_m\frac{(\sigma^{p+1}(t))^m}{m!}.\]
Since
\[(1+v)^x=\sum_{m=0}^\infty(x)_m\frac{v^m}{m!},\]
we can continue as
\[\sum_{m=0}^\infty(x)_m\frac{(\sigma^{p+1}(t))^m}{m!}=(1+\sigma^{p+1}(t))^x=(1+e^{\sigma^p(t)}-1)^x=(E^p(t))^x.\]
Definition \eqref{Belldef} and \eqref{mat2} gives the wanted identity.

To see the second one,
\[\frac{(\sigma^p(x))^m}{m!}=\frac{(e^{\sigma^{p-1}(x)}-1)^m}{m!}=\frac{1}{m!}\sum_{k=0}^m\binom{m}{k}(-1)^{m-k}[e^{\sigma^{p-1}(x)}]^k=\]
\[=\frac{1}{m!}\sum_{k=0}^m\binom{m}{k}(-1)^{m-k}[E^{p-1}(x)]^k.\]
Comparing the coefficients, the result follows from Proposition \ref{gfs}.
\hfill\qed

In his paper, Bell cited Dobi\'nski's result \cite{Dob}
\[(B_n^1=)\xi_n^{(1)}=\frac{1}{e}\sum_{k=0}^\infty\frac{k^n}{k!},\]
and wrote that ``there are similar (but more complicated) summations involving $\xi_n^{(m)}$, $m\ge1$". Without cited references with respect to the mentioned formulas, we do not know what would be these complicated sums, but we point out that beside Dobi\'nski's nice formula, the same and simple summation exist for higher order Bell numbers. What is more, Ces\`aro's integral representation formula \cite{Cesaro} can also be carried to the more general setting.

\begin{Proposition}\label{BellSum}For $p\ge 1$ we have
\[B_n^p(x)=\frac{1}{e^x}\sum_{k=0}^\infty\frac{x^kB_n^{p-1}(k)}{k!}.\]
\end{Proposition}

\textit{Proof.} The case $p=1$ is the result of Dobi\'nski. Let $p\ge1$. By the definition of $B_n^p(x)$,
\[\sum_{n=0}^\infty B_n^p(x)\frac{t^n}{n!}=(E^p(t))^x=e^{x\sigma^p(t)}=e^{x(e^{\sigma^{p-1}(t)}-1)}=\]
\[=\frac{1}{e^x}\sum_{k=0}^\infty\frac{x^k[e^{\sigma^{p-1}(t)}]^k}{k!}=\frac{1}{e^x}\sum_{k=0}^\infty\frac{x^k[E^{p-1}(t)]^k}{k!}.\]
\hfill\qed

And the integral representation of higher order Bell numbers is as follows.

\begin{Proposition}For $p\ge0$ we have
\[B_n^p(x)=\frac{2n!}{\pi}\im\int_0^\pi\exp(x\sigma^p(e^{it}))\sin(nt)dt.\]
\end{Proposition}

\textit{Proof.} The case $p=0$ comes from the orthogonality of sines on $[0,\pi]$ and the Taylor expansion of $\exp$ and the DeMoivre's formula. $p>0$ is the formula of Ces\`aro, so let $p>1$.
\begin{eqnarray*}
B_n^p(x)&=&\frac{1}{e^x}\sum_{k=0}^\infty\frac{x^kB_n^{p-1}(k)}{k!}=\\
&&\frac{2n!}{\pi}\sum_{k=0}^\infty\frac{x^k}{k!}\im\int_0^\pi\frac{\exp(k\sigma^{p-1}(e^{it}))}{e^x}\sin(nt)dt=\\
&&\frac{2n!}{\pi}\im\int_0^\pi\sum_{k=0}^\infty\frac{[x\exp(\sigma^{p-1}(e^{it}))]^k}{e^xk!}\sin(nt)dt=\\
&&\frac{2n!}{\pi}\im\int_0^\pi\exp(x(\exp(\sigma^{p-1}(e^{it}))-1))\sin(nt)dt=\\
&&\frac{2n!}{\pi}\im\int_0^\pi\exp(x\sigma^p(e^{it}))\sin(nt)dt.
\end{eqnarray*}
\hfill\qed

\section{Higher order Fubini numbers}

An other direction of generalizing old identities arises if we consider the ``ordered version" of Bell numbers.

The ordered Bell number (also known as Fubini number \cite{Gross,James,Tanny}) $F_n$, enumerates the ordered partitions of an $n$-set. Thus, we can define $F_n$'s as the values of the polynomials
\[\begin{pmatrix}F_0(x)\\F_1(x)\\\vdots\end{pmatrix}:=\mS\begin{pmatrix}0!x^0\\1!x^1\\\vdots\end{pmatrix}\]
at $x=1$. The generating function is \cite{Tanny,Sprugnoli}
\begin{equation}
\sum_{n=0}^\infty F_n(x)\frac{t^n}{n!}=\frac{1}{1-x(e^t-1)}.\label{Fubegf}
\end{equation}

Nice identities are valid for these numbers and polynomials \cite{Gross,James,Tanny,Wilf}:
\begin{equation}
F_n=\sum_{k=0}^\infty\frac{k^n}{2^{k+1}},\label{Fubinisum}
\end{equation}
\begin{equation}
F_n(x)=\sum_{k=0}^n A_{n,k}(x+1)^kx^{n-k},\label{FubEuler}
\end{equation}
where $A_{n,k}$ is the Eulerian number with parameters $n$ and $k$ and counts the permutations consisting $k$ rises on $n$ elements (see \cite{GKP}).

In what follows, we show that these relations remain true in the more general context.

\begin{Definition}\label{pFubdef}Let $p\ge1$, and let
\[\begin{pmatrix}F_0^p(x)\\F_1^p(x)\\\vdots\end{pmatrix}:=\mS^p\begin{pmatrix}0!x^0\\1!x^1\\\vdots\end{pmatrix}\]
We call these polynomials as higher order Fubini polynomials, and $F_n^p(1)=:F_n^p$ as higher order Fubini numbers.
\end{Definition}

It is very interesting that these numbers also appear in \cite{Blasetal}.

\begin{center}
\begin{tabular}{r|cccccccc}$F_n^p$&$n=1$&$n=2$&$n=3$&$n=4$&$n=5$&$n=6$&$n=7$&\\\hline
$p=1$&1&3&13&75&541&4683&47 293\\
$p=2$&1&4&23&175&1662&18 937&251 729\\
$p=3$&1&5&36&342&4048&57 437&950 512\\
$p=4$&1&6&52&594&8444&143 783&2 854 261\\
$p=5$&1&7&71&949&15 775&313 920&7 279 795\\
\end{tabular}
\end{center}

\begin{Proposition}With the above definition, the Fubini numbers of higher order can be represented as
\[F_n^p=\sum_{k=0}^\infty\frac{B_n^p(k)}{2^{k+1}}.\]
\end{Proposition}

\textit{Proof.} Rewrite Proposition \ref{PropBellSt}. as
\[\sum_{m=0}^nS_{n,m}^px^m=B_n^p(x).\]
If we apply Proposition \ref{BellSum}, multiply both sides with $1/e^x$ and integrate on $[0,+\infty[$ with respect to $x$, we get
\[\sum_{m=0}^nS_{n,m}^p\int_0^\infty\frac{x^m}{e^x}=\int_0^\infty\frac{1}{e^{2x}}\sum_{k=0}^\infty\frac{x^kB_n^{p-1}(k)}{k!}.\]
The integral on the left hand side is $m!$, and
\[\int_0^\infty\frac{x^k}{e^{2x}}=\frac{k!}{2^{k+1}},\]
hence the result follows.
\hfill\qed

It is not hard to deduce the exponential generating functions for these numbers.

\begin{Proposition}We have for all $p\ge1$
\[\sum_{n=0}^\infty F_n^p(x)\frac{t^n}{n!}=\frac{1}{1-x\sigma^p(t)}.\]
\end{Proposition}

\textit{Proof.} Case $p=1$ is known (see \eqref{Fubegf} and the references before the formula), so let $p>1$ be fixed. We define a linear operator $L$ on $\RR[[x]]$ as
\[L((x)_m):=m!x^m.\]
(This defines $L$ entirely, because the sequence $((x)_m)_{m=0}^\infty$ is a base in $\RR[[x]]$.)
Then, we apply $L$ to \eqref{BellStsum}:
\[L(B_n^{p-1}(x))=\sum_{m=0}^n S_{n,m}^pL((x)_m)=\sum_{m=0}^n S_{n,m}^pm!x^m=F_n^p(x).\]
This gives
\[\sum_{n=0}^\infty F_n^p(x)\frac{t^n}{n!}=\sum_{n=0}^\infty L(B_n^{p-1}(x))\frac{t^n}{n!}=L\left(\sum_{n=0}^\infty B_n^{p-1}(x)\frac{t^n}{n!}\right)=L\left(E^{p-1}(t))^x\right).\]
Now we introduce the temporary variable $v:=E^{p-1}(t)-1$ ($t$ can be considered as a fixed real number). Whence
\[L\left(E^{p-1}(t))^x\right)=L((1+v)^x)=L\left(\sum_{n=0}^\infty(x)_n\frac{v^n}{n!}\right)=\sum_{n=0}^\infty L((x)_n)\frac{v^n}{n!}=\sum_{n=0}^\infty\frac{(xv)^n}{n!}=\]
\[\frac{1}{1-xv}=\frac{1}{1-x(E^{p-1}(t)-1)}=\frac{1}{1-x\sigma^p(t)}.\]
This proves the proposition.
\hfill\qed

\section{Higher order Eulerian numbers}

Now we show that how can one construct higher order Eulerian numbers. Let us recall \eqref{FubEuler}, and write it in the form
\[\begin{pmatrix}\frac{F_0(x)}{x^0}\\\frac{F_1(x)}{x^1}\\\vdots\end{pmatrix}:=\mathcal{A}\begin{pmatrix}\left(\frac{x+1}{x}\right)^0\\\left(\frac{x+1}{x}\right)^1\\\vdots\end{pmatrix}.\]
Here the $(n,k)$-th entry of $\mathcal{A}$ is $A_{n,k}$, the Eulerian number with parameters $n$ and $k$. According to Definition \ref{pFubdef} and identity \eqref{BellStsum}, it seems to be fruitful to give a similar

\begin{Definition}Let $\mathcal{A}^{(p)}$ be the (unique) matrix $(p\ge1)$, for which
\[\begin{pmatrix}\frac{F_0^p(x)}{x^0}\\\frac{F_1^p(x)}{x^1}\\\vdots\end{pmatrix}:=\mathcal{A}^{(p)}\begin{pmatrix}\left(\frac{x+1}{x}\right)^0\\\left(\frac{x+1}{x}\right)^1\\\vdots\end{pmatrix}.\]
In addition, we introduce the polynomials
\[\begin{pmatrix}A_0^p(x)\\A_1^p(x)\\\vdots\end{pmatrix}:=\mathcal{A}^{(p)}\begin{pmatrix}x^0\\x^1\\\vdots\end{pmatrix}.\]
\end{Definition}

We shall point out that this definition is the proper one, since it gives back all the existing identities connecting Fubini, Stirling and Eulerian numbers and polynomials (Proposition \ref{PropEulerian}.).

It is important to notice that $\mathcal{A}^{(p)}$ is not equal to the $p$-th power of $\mathcal{A}$. This is why we use the brackets.

Because the identity \eqref{FubEuler} was generalized on an other direction in a former paper of the author in \cite{Mezo}, all the proofs of the next Proposition can be deduced verbatim. These identities are straight generalizations of the ordinary case $(p=1)$.

\begin{Proposition}\label{PropEulerian}The next identities are valid
\begin{eqnarray*}
m!S_{n,m}^p&=&\sum_{k=0}^nA_{n,k}^p\binom{k}{n-m},\\
A_{n,m}^p&=&\sum_{k=0}^nk!S_{n,k}^p\binom{n-k}{m}(-1)^{n-k-m},\\
F_n^p(x)&=&\sum_{k=0}^nA_{n,k}^p(x+1)^kx^{n-k},\\
F_n^p(x)&=&x^nA_n^p\left(\frac{x+1}{x}\right).
\end{eqnarray*}
\end{Proposition}

\begin{Corollary}The exponential generating function for the polynomial $A_n^p(x)$ is
\[\sum_{n=0}^\infty A_n^p(x)\frac{t^n}{n!}=\frac{1}{1-\frac{1}{x-1}\sigma^p((x-1)t)}=\frac{x-1}{x-E^{p-1}((x-1)t)}.\]
\end{Corollary}

\section{Application to the eigensequence of the Stirling transform}

In this section we present a formula for the eigensequence of the Stirling transform involving $S^2_{n,m}$. First, let us introduce the necessary notions. Let a sequence $a_n$ be given. Then the Stirling transform $b_n$ of $a_n$ is defined by
\[b_n=\sum_{k=0}^na_kS_{n,k}\quad(n=0,1,\dots).\]
We consider the eigensequence of the Stirling transform \cite{BerSloane}, that is, $C_n$ is defined as
\[C_{n+1}=\sum_{k=0}^nC_kS_{n,k}.\]
(The sequence begins with $0,1,1,2,6,26,152,1144,\dots$.) With our notation,
\[\mS\begin{pmatrix}C_0\\C_1\\\vdots\end{pmatrix}=\begin{pmatrix}C_1\\C_2\\\vdots\end{pmatrix}\]
It is known \cite{BerSloane} that the exponential generating function of $C_n$ satisfies the equation
\[f'(x)=f(e^x-1)+1,\]
or, with the notations have introduced,
\[f'(x)=f(\sigma^1(x))+1.\]
Then it can be readily checked that
\begin{equation}
f''(x)=e^x\left[f(\sigma^2(x))+1\right].\label{Angf}
\end{equation}
We see that
\[f(\sigma^2(x))=\sum_{k=0}^\infty C_k\frac{(\sigma^2(x))^k}{k!}=\sum_{k=0}^\infty C_k\sum_{n=0}^\infty S^2_{n,k}\frac{x^n}{n!}=\sum_{n=0}^\infty\frac{x^n}{n!}\left(\sum_{k=0}^\infty C_kS^2_{n,k}\right).\]
The coefficients of $x^n/n!$ are the elements of the vector
\[\mS^2\begin{pmatrix}C_0\\C_1\\\vdots\end{pmatrix}=:\begin{pmatrix}C_0'\\C_1'\\\vdots\end{pmatrix}\]
(truncated after $n$).
The final consideration is that
\[e^x\left[f(\sigma^2(x))+1\right]=e^x\sum_{n=0}^\infty C_n'\frac{x^n}{n!}+\sum_{n=0}^\infty\frac{x^n}{n!}=\sum_{n=0}^\infty\frac{x^n}{n!}\left[1+\sum_{k=0}^n\binom{n}{k}C'_k\right].\]
Hence \eqref{Angf} implies the identity
\[C_{n+2}=1+\sum_{k=0}^n\binom{n}{k}C'_k,\]
or, rephrased applying the matrix formalism,
\[\begin{pmatrix}C_2\\C_3\\\vdots\end{pmatrix}=\textbf{1}+\mathcal{B}\begin{pmatrix}C_0'\\C_1'\\\vdots\end{pmatrix}=\textbf{1}+\mathcal{B}\mS^2\begin{pmatrix}C_0\\C_1\\\vdots\end{pmatrix}.\]
Here
\[\textbf{1}=\begin{pmatrix}1\\1\\\vdots\end{pmatrix},\quad\mathcal{B}=\begin{pmatrix}\binom{0}{0}&0&0&0&0&\cdots\\\binom{1}{0}&\binom{1}{1}&0&0&0&\cdots\\\binom{2}{0}&\binom{2}{1}&\binom{2}{2}&0&0&\cdots\\\vdots&\vdots&\ddots&\vdots&\vdots\\\binom{n}{0}&\binom{n}{1}&\cdots&\binom{n}{n}&0&\cdots\\\vdots\end{pmatrix}.\]

In addition, we prove the next
\begin{Proposition}We have
\[C_{n+2}=\sum_{k=0}^{n}S_{n,k}(C_{k+1}+kC_k),\]
that is,
\[\begin{pmatrix}C_2\\C_3\\\vdots\end{pmatrix}=\mS\begin{pmatrix}C_1\\C_2\\\vdots\end{pmatrix}+\mS\begin{pmatrix}0\cdot C_0\\1\cdot C_1\\\vdots\end{pmatrix}.\]
\end{Proposition}

\textit{Proof.} Let us introduce the operator $\mathcal{L}$ as
\[\mathcal{L}\begin{pmatrix}x^0\\x^1\\\vdots\end{pmatrix}=\begin{pmatrix}C_0\\C_1\\\vdots\end{pmatrix}.\]
Then
\[\mathcal{L}\mS\begin{pmatrix}x^0\\x^1\\\vdots\end{pmatrix}=\mS\mathcal{L}\begin{pmatrix}x^0\\x^1\\\vdots\end{pmatrix}=\begin{pmatrix}C_1\\C_2\\\vdots\end{pmatrix}.\]
Recalling \eqref{Stoxpow}, we get
\[\mathcal{L}\begin{pmatrix}B_0(x)\\B_1(x)\\\vdots\end{pmatrix}=\begin{pmatrix}C_1\\C_2\\\vdots\end{pmatrix}.\]
So, by \cite[p. 114]{Aigner},
\[C_{n+2}=\mathcal{L}B_{n+1}(x)=\mathcal{L}(xB_n(x)+xB_n'(x))=\mathcal{L}(xB_n(x))+\mathcal{L}(xB_n'(x))=\]
\[\mathcal{L}\left(\sum_{n=0}^nS_{n,k}x^{k+1}\right)+\mathcal{L}\left(\sum_{n=0}^nS_{n,k}kx^k\right)=\sum_{n=0}^nS_{n,k}\mathcal{L}(x^{k+1})+\sum_{n=0}^nS_{n,k}k\mathcal{L}(x^k),\]
and the result follows.
\hfill\qed

\end{document}